\documentclass[]{opt2020} 

\usepackage{booktabs,bbm,amsmath}
\usepackage{xcolor,bm,amsmath,graphicx,comment,bigints}
\usepackage{algorithm,algorithmic}

\graphicspath{{images/}}

\title[A Homotopy Algorithm for Optimal Transport]{A Homotopy Algorithm for Optimal Transport}

\optauthor{\Name{Roozbeh Yousefzadeh} \Email{roozbeh.yousefzadeh@yale.edu}\\
  \addr Yale University and VA Connecticut Healthcare System, New Haven, CT, USA}

\begin{document}

\maketitle

\begin{abstract}%
The optimal transport problem has many applications in machine learning, physics, biology, economics, etc. Although its goal is very clear and mathematically well-defined, finding its optimal solution can be challenging for large datasets in high-dimensional space. Here, we propose a homotopy algorithm that first transforms the problem into an easy form, by changing the target distribution. It then transforms the problem back to the original form through a series of iterations, tracing a path of solutions until it finds the optimal solution for the original problem. We define the homotopy path as a subspace rotation based on the orthogonal Procrustes problem, and then we discretize the homotopy path using eigenvalue decomposition of the rotation matrix. Our goal is to provide an algorithm with complexity bound $\mathcal{O}(n^2 \log(n))$, faster than the existing methods in the literature.
\end{abstract}


\section{Introduction} \label{sec:intro}

The optimal transport problem has been the subject of many studies in recent years as it has several applications in machine learning and other fields, such as physics \citep{buttazzo2012optimal}, biology \citep{schiebinger2019optimal}, economics \citep{galichon2018optimal}, etc. Study of this problem has, in fact, a long history, going back to the 18th century \citep{villani2008optimal}. Although there are well-known algorithms for solving this problem, those algorithms are not practical or fast enough when we deal with modern machine learning datasets with large number of samples and many features \citep{peyre2019computational}. In fact, the best known complexity bound for this problem is $\mathcal{O}(n^3)$, for distributions with $n$ samples \cite{pele2009fast}. In pursuit of faster algorithms, a regularized (relaxed) version of this problem was proposed by \cite{cuturi2013sinkhorn} and there has been many advances in solving the regularized formulation, e.g., \cite{lin2020projection,lin2019efficient}. Here, we design the outline of a homotopy algorithm for the original formulation, suggesting that the optimal solution may be computed with complexity bound $\mathcal{O}(n^2 \log(n))$. Our algorithm exploits the structure of the data and starts by solving a related problem with known optimal solution. By following a continuous path, it then transforms the related problem back to its original form until it finds the solution for the original problem.

\section{Formulation} \label{sec:formulation}

We would like to find a 1-to-1 mapping from one set of points $\mathcal{X}$ (a discrete distribution) to another set of points $\mathcal{Y}$, such that the sum of distances between the matched points is minimal. Each set has $n$ discrete points in $\mathbb{R}^d$. In the general form of the problem, the distance between the points can be measured using different metrics \citep{gramfort2015fast}. In our formulation, we use the Euclidean distance which is common in the optimal transport literature and is used for analyzing the complexity of most recent algorithms  \citep{lin2020projection,lin2019efficient,guo2020fast}. 
In this setting, each point set can be considered a discrete probability distribution. Using a transportation plan $T$, the goal is to basically transport the discrete distribution $\mathcal{X}$ to the discrete distribution $\mathcal{Y}$.

We organize the information about the discrete distributions $\mathcal{X}$ and $\mathcal{Y}$ in matrix form, where columns represent the $n$ points (i.e., samples) and rows represent the $d$ dimensions (i.e, features). 
$T$ is a binary $n\times n$ matrix that determines the pairwise matching between the points in $\mathcal{X}$ and $\mathcal{Y}$. Since the mapping is 1-to-1, the sum of each row and each column of $T$ should be equal to 1. 
Our optimization problem is as following:
\vspace{-.2cm}
\begin{equation} \label{eq:obj}
    \min_{T} \| \mathcal{X} - \mathcal{Y} T \|_F,
\end{equation}
subject to:
\begin{equation} \label{eq:const1}
    T \mathbbm{1}_{n,1} = \mathbbm{1}_{n,1},
\end{equation}
\vspace{-.5cm}
\begin{equation} \label{eq:const2}
    \mathbbm{1}_{1,n} T = \mathbbm{1}_{1,n},
\end{equation}
\vspace{-.5cm}
\begin{equation} \label{eq:const3}
    T \in \{0,1\}.
\end{equation}

\vspace{-.5cm}


\section{Our homotopy approach}

We now design our homotopy algorithm to solve the above formulation. In numerical optimization, homotopy algorithms are usually used to solve hard optimization problems \citep{nocedal2006numerical}. In such approach, a homotopy transformation is designed to overcome the specific challenges of the hard problem. 
Recent examples of developing homotopy algorithms include \citep{dunlavy2005homotopy,mobahi2015link,anandkumar2017homotopy,yousefzadeh2020deep}. 
In our approach, we transform the target distribution $\mathcal{Y}$ into a distribution $\mathcal{Y}_h$ such that the solution to our optimization problem (Equations~\eqref{eq:obj}-\eqref{eq:const3}) is unique and easily computable. We then gradually transform the $\mathcal{Y}_h$ back to its original form, $\mathcal{Y}$, through a series of iterations.

\subsection{Investigating the formulation}

Let's start by investigating our objective function and find a plausible way to transform it into a state where the optimal solution is known. From linear algebra, it is easy to show that

\begin{equation} \label{eq:bound}
    \| \mathcal{X} - \mathcal{Y} T \|_F^2 = \| \mathcal{X} \|_F^2 + \| \mathcal{Y} \|_F^2 - 2 tr(T^T \mathcal{Y}^T \mathcal{X}).
\end{equation}
where $tr(.)$ is the trace operator. 
Our objective function is, therefore, equivalent to

\begin{equation} \label{eq:objmax}
    \max_T tr(T^T \mathcal{Y}^T \mathcal{X}),
\end{equation}
subject to constraints~\eqref{eq:const1}-\eqref{eq:const3}. This is finding a permutation matrix that maximizes the trace of $\mathcal{Y}^T\mathcal{X}$. 

\subsection{Homotopy transformation}

To design a homotopy algorithm, we have to first design a transformation that makes our optimization problem easier to solve. Besides that, we also want the transformed problem to be related to the original problem, so that the solution of transformed problem can lead us to the solution of original problem.

Based on the above considerations, we replace the $T$ with $Q$, and instead of constraints~\eqref{eq:const1}-\eqref{eq:const3}, we only require the $Q$ to be orthogonal, i.e., $Q^T Q = I$. With such relaxation, we can analytically find its unique optimal solution, $Q^*$. In fact, this relaxation would provide a lower bound to the solution of our minimization problem in Equation~\eqref{eq:obj}, which would translate to an upper bound on the solution of our maximization problem in Equation~\eqref{eq:objmax}. Note that $T$ is also orthogonal, so in our transformation, we have preserved that property for $Q$. As we will see later, orthogonality of $Q$ will be helpful in discretizing the homotopy path.

The way to compute the $Q$ is to compute the SVD of $\mathcal{Y}^T \mathcal{X} = U \Sigma V^T$. We can then compute the optimal solution as $Q^* = U V^T$, as proved in \cite[Section 6.4]{golub2012matrix}. 
Now, if we replace our $\mathcal{Y}$ with $\mathcal{Y}_h = \mathcal{Y} Q^*$ and set $T = \mathcal{I}$, we have actually optimized our problem for distributions $\mathcal{X}$ and $\mathcal{Y}_h$. This defines our homotopy transformation: transforming $\mathcal{Y}$ to $\mathcal{Y}_h$. 
This transformation is merely a step to make the problem tractable, so we do not need to precisely compute the $Q^*$. Instead, we can use a randomized SVD algorithm to approximate the $Q^*$. The complexity of computing the $Q^*$ is therefore, $\mathcal{O}(n^2 \log(k))$, where $k$ is the number of singular values used in the approximation \citep{halko2011finding}.

\subsection{Homotopy path}

The next step is to follow a homotopy path to gradually transform the $\mathcal{Y}_h$ back to its original form, $\mathcal{Y}$. At each iteration $i$, $\mathcal{Y}_i$ moves closer towards $\mathcal{Y}$, and we optimize the $T^i$ for $\mathcal{Y}_i$, using the $T^{i+1}$ as the starting point. At the end, we obtain a solution for transporting $\mathcal{X}$ to~$\mathcal{Y}$.

One can follow different paths to transition $\mathcal{Y}_h$ back to $\mathcal{Y}$. Here, we choose a homotopy path that corresponds to the transformation used to obtain the $\mathcal{Y}_h$ in the first place. The transformation performed with $Q^*$ is a subspace rotation, known in linear algebra as the "orthogonal Procrustes problem". Since our homotopy is a rotation, we devise a path to gradually reverse the rotation.

We note that $Q^*$ is orthogonal and non-singular. Therefore, using a blocked Schur algorithm \citep{deadman2012blocked}, we can compute its square root, ${Q^*}^{1/2}$, such that it is a function of $Q^*$. Such square root is guaranteed to exist for $Q^*$ and it would also be orthogonal \citep{higham1987computing,higham2005functions}. This breaks our overall rotation into two steps. 
By computing the square root of ${Q^*}^{1/2}$, we can then break each of the above two steps into two finer steps, and so on. This way, we can discretize our continuous homotopy path. As an example, if we discretize 3 times, our incremental rotation matrix will be ${Q^*}^{1/8}$. This corresponds to $2^3$ steps to move from $\mathcal{Y}_h$ to $\mathcal{Y}$.

Since the square roots are a function of $Q^*$, our homotopy path is also a function of $Q^*$, and our discretization points will be on a continuous homotopy path. Moreover, since $Q^*$ (and its square roots) is orthogonal, its condition number is 1 which is desirable for numerical stability.

We note that following this procedure, the number of steps will be a power of 2. 
The best step size would be such that the solution at each step is relatively close to the solution from the previous step, so that we can easily move from one solution to another. Making the step sizes very fine can make the process inefficient, as the solution might stay the same between very fine steps.

Regarding the computational complexity, the number of homotopy steps, $h$, is a constant, e.g., 2, 4, 8, or 16. Hence, it does not add to the overall complexity of our algorithm. The complexity of computing the square roots of $Q^*$ is equivalent to the complexity of eigenvalue decomposition, $\mathcal{O}(n^2 \log(n))$.

\subsection{Updating the solution along the path}

At each iteration on the homotopy path, we aim to improve the solution from the previous step. For that, we use an optimization module
\vspace{-.2cm}

\begin{equation}
    T = \mathcal{F}(\mathcal{X},\mathcal{Y},T_{initial}),
\end{equation}
which takes as input the two discrete distributions and a $T_{initial}$ as the initial solution for optimization variables. It then returns $T$ as the optimal solution. 

Since $T$ is a permutation matrix, $\mathcal{F}$ should be a combinatorial algorithm. There are off-the-shelf exact and approximation algorithms for $\mathcal{F}$. Since we are following a homotopy path, we expect the optimal solution at each iteration to be relatively similar to the previous iteration. Therefore, we can rely on efficient algorithms that only aim to improve the solution by making local modifications to it. For example, \cite{duan2014linear} and \cite{altschuler2017near} perform local updates in near-linear time with approximation guarantees.

Comprehensive options for $\mathcal{F}$ aim to optimize the $T$ globally, e.g. \cite{duff2001algorithms}. We use them to verify our results and to ensure that our local updates are adequate. Algorithm~\ref{alg:homotopy} formalizes these procedures.

\begin{algorithm}[H]
\caption{Homotopy Algorithm for Optimal Transport}
\label{alg:main}
\begin{small}
\textbf{Inputs}: Discrete distribution $\mathcal{X}$, target distribution $\mathcal{Y}$, number of homotopy steps $h$ (power of 2)\\
\textbf{Outputs}: Transport plan $T$, total transport cost $\kappa$
\begin{algorithmic}[1] 
\label{alg:homotopy}
\STATE $G = greedysort(\mathcal{X},\mathcal{Y})$
\STATE $\mathcal{Y}_g = \mathcal{Y} G$
\STATE $U \Sigma V^T = svd(\mathcal{Y}_g^T \mathcal{X})$
\STATE $T_h = \mathcal{I}_{n \times n}$
\STATE $T_{\delta} = \sqrt[h]{U V^T}$
\FOR{$i=h-1$ to $0$}
    \STATE $\mathcal{Y}_i = real( \mathcal{Y}_g (T_{\delta})^{i} )$
    \STATE $T_i = \mathcal{F}(\mathcal{X},\mathcal{Y}_i,T_{i+1})$
\ENDFOR
\STATE $T = G T_0$
\STATE $\kappa = \| \mathcal{X} - \mathcal{Y} T \|_F$
\STATE \textbf{return} $T$ and $\kappa$
\end{algorithmic}
\end{small}
\end{algorithm}

$greedysort(.)$ in Algorithm~\ref{alg:homotopy} basically performs a pass over the samples in $\mathcal{X}$. For each sample in $\mathcal{X}$, it assigns the closest unassigned sample in $\mathcal{Y}$, until it has a 1-to-1 matching between the points in $\mathcal{X}$ and $\mathcal{Y}$. This has complexity $\mathcal{O}(n)$ and it is a quick and efficient attempt to initialize the transportation plan. 

We note that the homotopy transformation has overall complexity of $\mathcal{O}(n^2log(n))$. The process of updating the optimal solutions along the homotopy path can also be performed in $\mathcal{O}(n^2log(n))$ using approximation algorithms for $\mathcal{F}$, leading to overall complexity of $\mathcal{O}(n^2log(n))$ for Algorithm~\ref{alg:homotopy}.

\vspace{-.5cm}

\section{Toy example} \label{sec:toy}

\begin{figure}[h]
  \centering
   \includegraphics[width=0.25\linewidth]{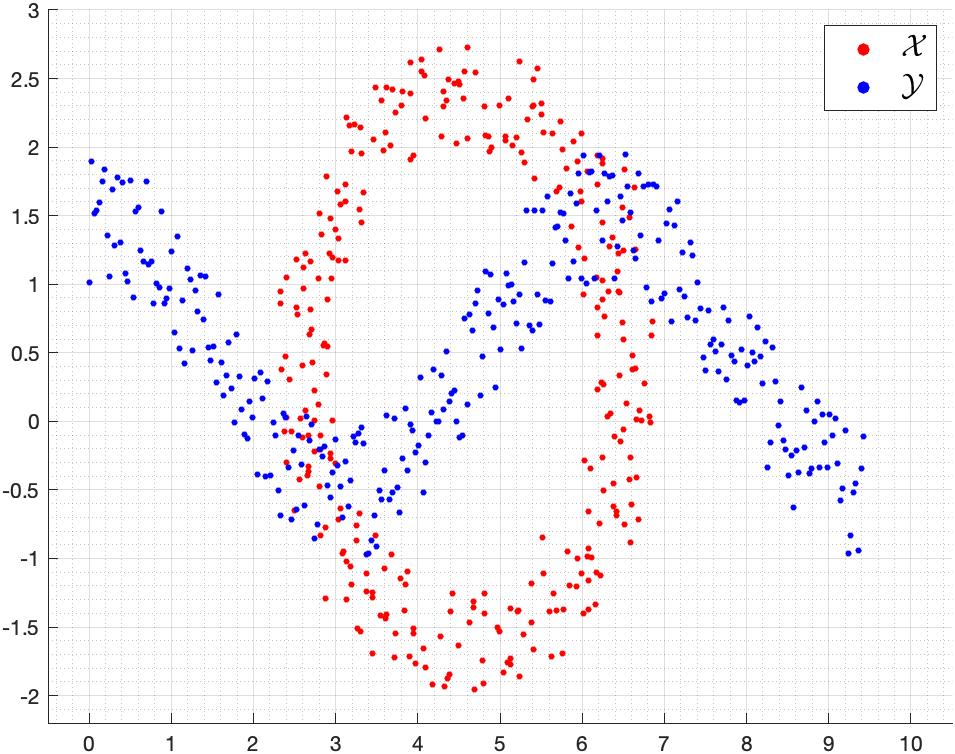}
   \caption{Demonstration of two distributions for our toy example: $\mathcal{X}$ in red and $\mathcal{Y}$ in blue.}
  \label{fig:toy_data}
\end{figure}

We start by a small dataset of 300 data points in 2-dimensional space, in order to demonstrate how our algorithm works. Figure~\ref{fig:toy_data} demonstrates our two discrete distributions: $\mathcal{X}$ in red and $\mathcal{Y}$ in blue.

\vspace{-.2cm}

\begin{figure}[h]
  \centering
   \includegraphics[width=0.25\linewidth]{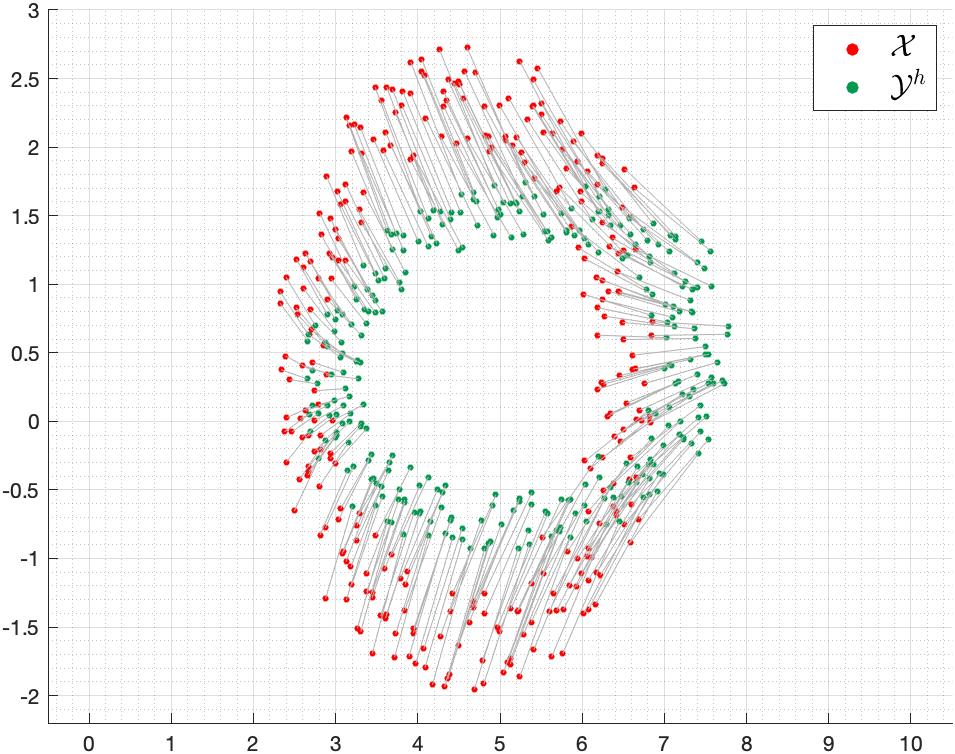}
   \caption{Transformed distribution $\mathcal{Y}_h$ shown in green.
   }
  \label{fig:toy_transformed}
\end{figure}

Our algorithm transforms the $\mathcal{Y}$ distribution into $\mathcal{Y}_h$ which is shown in green in Figure~\ref{fig:toy_transformed}. The gray lines in this figure show the optimal matching between the $\mathcal{X}$ and $\mathcal{Y}_h$ leading to total transport cost $\kappa = 15.47$. This can be considered a lower bound on the cost.

The next step is to gradually transform the $\mathcal{Y}_h$ back into its original form $\mathcal{Y}$. We use $h=4$ as the number of steps. Figure~\ref{fig:toy_walkback} shows the target distributions ($\mathcal{Y}_i$) at each of the intermediate states.

\begin{figure}[h]
\floatconts
  {fig:toy_walkback}
  {    \vspace{-.5cm}
\caption{Gradual transformation of $\mathcal{Y}_h$ back into its original form, $\mathcal{Y}$. Note that \textbf{(a)} depicts $\mathcal{Y}_h$ as in Figure~\ref{fig:toy_transformed}, and \textbf{(e)} depicts $\mathcal{Y}$ as in Figure~\ref{fig:toy_data}.}}
  {%
    \subfigure[]{\label{fig:circle}%
      \includegraphics[width=0.2\linewidth]{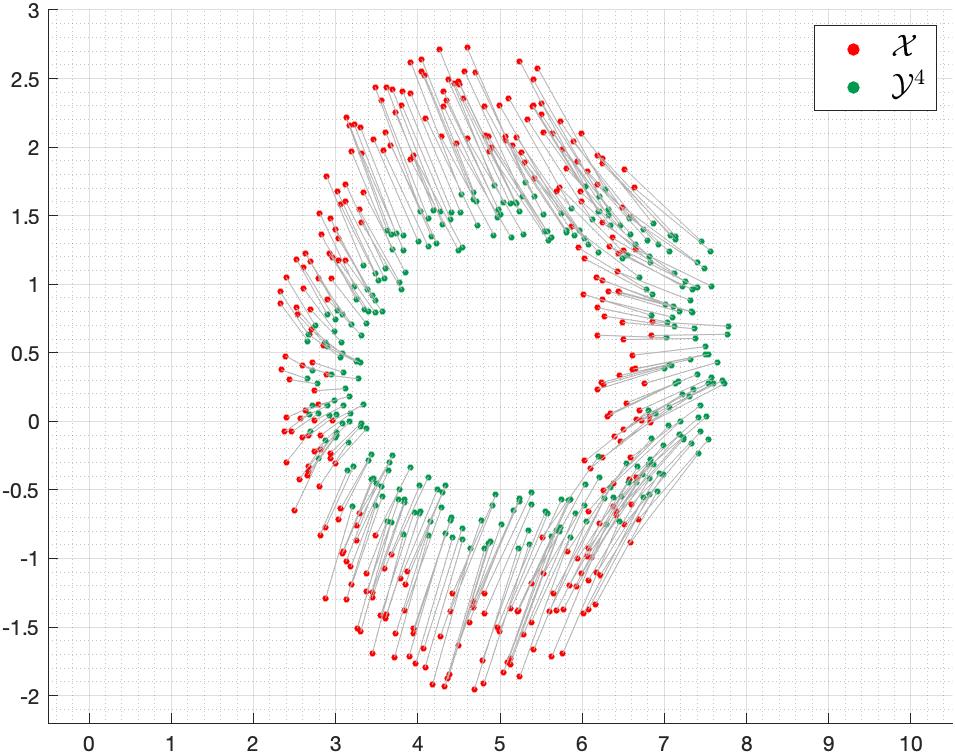}}
    \subfigure[]{\label{fig:square}%
      \includegraphics[width=0.2\linewidth]{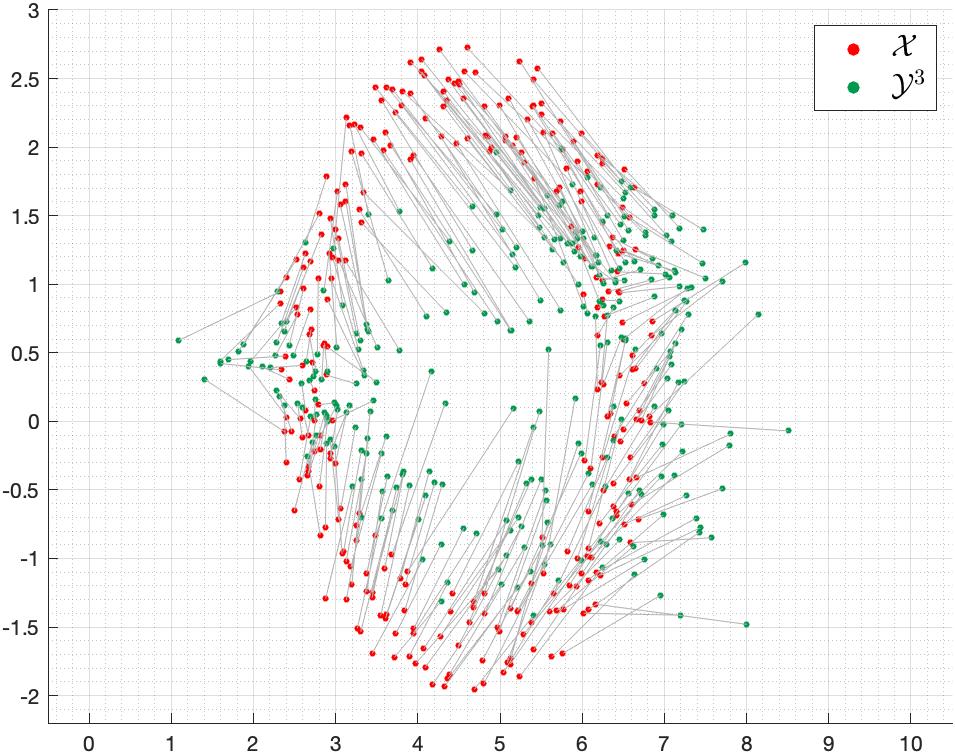}}
    \subfigure[]{\label{fig:square}%
      \includegraphics[width=0.2\linewidth]{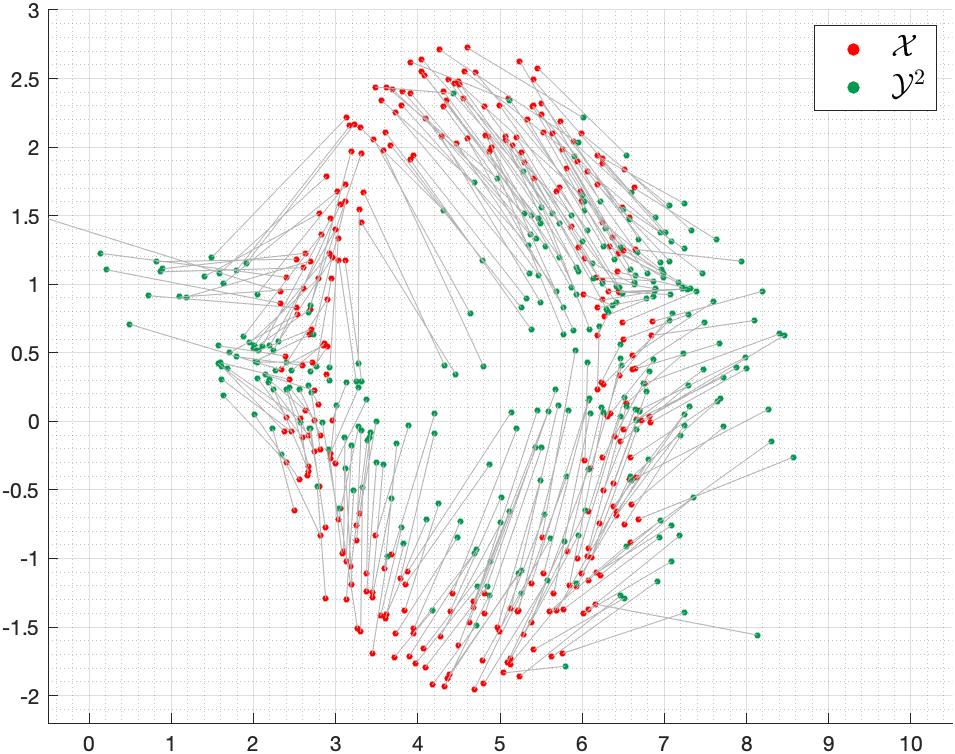}} \\
    \subfigure[]{\label{fig:square}%
      \includegraphics[width=0.2\linewidth]{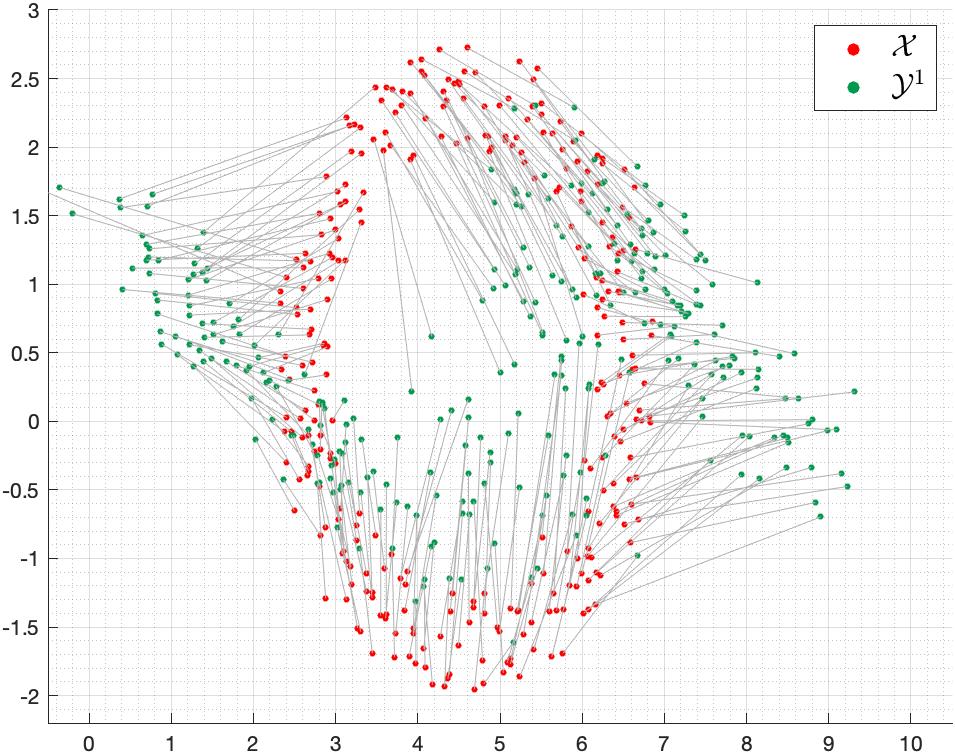}}
    \subfigure[]{\label{fig:square}%
      \includegraphics[width=0.2\linewidth]{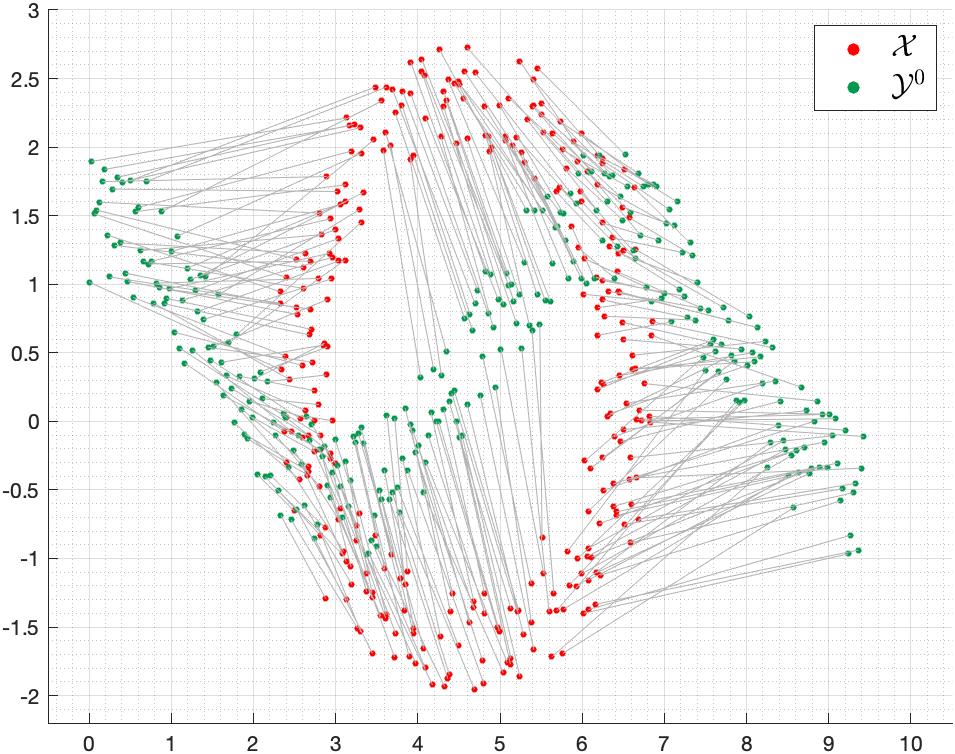}}}
\end{figure}


\begin{table}[h]
\caption{Results of Algorithm~\ref{alg:homotopy} along the homotopy path.}
\label{tbl:toy_results}
\vskip 0.15in
\begin{center}
\begin{footnotesize}
\begin{sc}
\begin{tabular}{c c c c}
\toprule
Iteration $(i)$  & $\kappa(\mathcal{X},\mathcal{Y}_i,T_{i+1})$  & $\kappa(\mathcal{X},\mathcal{Y}_i,T_{i})$  \\
\midrule
4  &  - & 15.47 \\
3  &  22.01 & 18.06 \\
2  &  26.10 & 20.86 \\
1  &  27.65 & 24.50 \\
0  &  31.23 & 29.06 \\
\bottomrule
\end{tabular}
\end{sc}
\end{footnotesize}
\end{center}
\end{table}

Table~\ref{tbl:toy_results} shows the results of optimizing the problem along the homotopy path. At each iteration, we update the $T_{i}$ for the intermediate distribution $\mathcal{Y}_i$, using a linear-time algorithm \citep{duan2014linear}. Using a more comprehensive algorithm (e.g., algorithm by \cite{duff2001algorithms} with cubic running time) yields the same results as when we use the linear-time algorithm, because we initiate the homotopy path with an optimal solution, and at each homotopy step, the solution from previous iteration is nearly optimal, so we only have to make small adjustments to it.


The third column in Table~\ref{tbl:toy_results} shows the total cost corresponding to the combination of $\mathcal{Y}_i$ and $T_{i+1}$ at each step. This basically means using the solution from the previous step for the current step. Clearly, the improvement achieved by $\mathcal{F}$ at each step is relatively small, because the solution from the previous step is almost optimal and considerably close to the optimal solution for the current step. Note that when we use a random permutation matrix for $T$, the total cost is above 60 on average.

When we apply the cubic running time algorithm by \cite{duff2001algorithms} directly to $\mathcal{X}$ and $\mathcal{Y}$, the optimal solution is the same as the optimal solution obtained by our homotopy algorithm. Also, at each iteration along the homotopy path, the cubic running time algorithm finds the same solution as the linear algorithm. Next, we investigate the effect of homotopy path discretization for this example (i.e., the number of homotopy steps).

\paragraph{Effect of the homotopy path discretization.}

To investigate the effect of the number of steps, $h$, we repeat this experiment with $h=2$ and $h = 8$. The results are shown in Tables~\ref{tbl:toy_results_h2} and~\ref{tbl:toy_results_h8}.

\begin{table}[h]
\caption{Results of Algorithm~\ref{alg:homotopy} along the homotopy path with $h=2$.}
\label{tbl:toy_results_h2}
\vskip 0.15in
\begin{center}
\begin{small}
\begin{sc}
\begin{tabular}{c c c c}
\toprule
Iteration $(i)$  & $\kappa(\mathcal{X},\mathcal{Y}_i,T_{i+1})$  & $\kappa(\mathcal{X},\mathcal{Y}_i,T_{i})$  \\
\midrule
2  &  - & 15.47 \\
1  &  33.63 & 20.86 \\
0  &  37.59 & 29.06 \\
\bottomrule
\end{tabular}
\end{sc}
\end{small}
\end{center}
\end{table}

\begin{table}[h]
\caption{Results of Algorithm~\ref{alg:homotopy} along the homotopy path with $h=8$.}
\label{tbl:toy_results_h8}
\vskip 0.15in
\begin{center}
\begin{small}
\begin{sc}
\begin{tabular}{c c c c}
\toprule
Iteration $(i)$  & $\kappa(\mathcal{X},\mathcal{Y}_i,T_{i+1})$  & $\kappa(\mathcal{X},\mathcal{Y}_i,T_{i})$  \\
\midrule
8  &  - & 15.47 \\
7  &  17.38 & 16.29 \\
6  &  19.30 & 18.06 \\
5  &  21.16 & 19.50 \\
4  &  22.10 & 20.86 \\
3  &  23.43 & 22.48 \\
2  &  25.32 & 24.50 \\
1  &  27.40 & 26.85 \\
0  &  29.72 & 29.06 \\
\bottomrule
\end{tabular}
\end{sc}
\end{small}
\end{center}
\end{table}

\begin{figure}[H]
  \centering
   \includegraphics[width=0.7\linewidth]{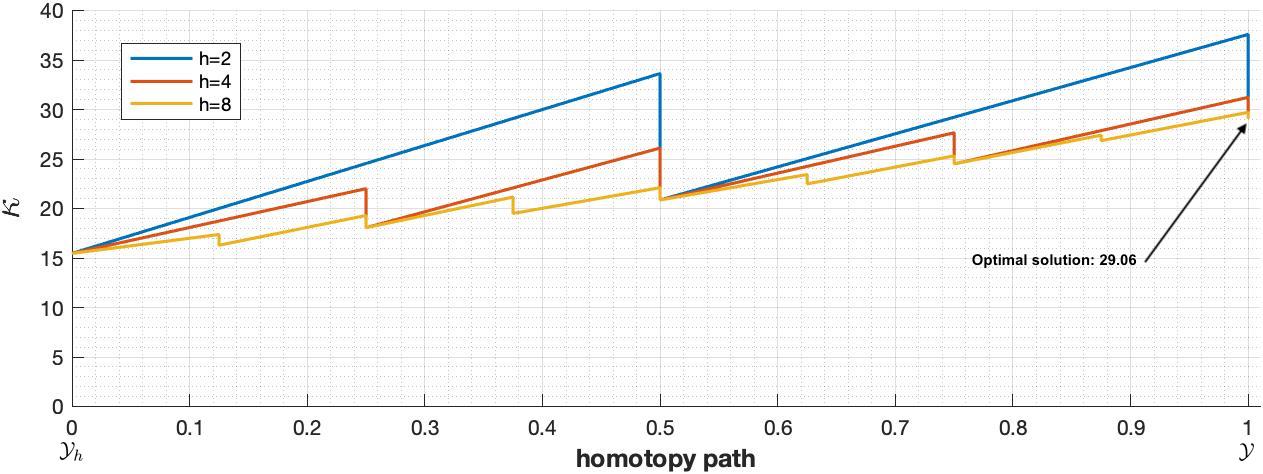}
   \caption{The effect of fineness of discretization of the homotopy path. The horizontal axis shows the path, which is normalized between 0 ($\mathcal{Y}_h$) and 1 ($\mathcal{Y}$). The vertical axis shows the value of $\kappa$ for the solutions we trace along the path. Note that we directly connect the discretization points using lines. Each vertical drop corresponds to an improvement achieved by $\mathcal{F}$.}
  \label{fig:toy_pathsteps}
\end{figure}

Figure~\ref{fig:toy_pathsteps} demonstrates the total cost associated with the solutions we trace on the homotopy path and how they are affected by the choice of $h$.

The horizontal axis in this Figure corresponds to the homotopy path, normalized between 0 ($\mathcal{Y}_h$) and 1 ($\mathcal{Y}$). The vertical axis shows the value of $\kappa$ for the solutions we trace along the path. 

As an example, let's walk along the path for the case of $h=2$, the blue line in Figure~\ref{fig:toy_pathsteps}. We start from 0 (far left of horizontal axis) on the homotopy path where the target distribution is $\mathcal{Y}_{2}$, the optimal solution is $T_2 = \mathcal{I}$, and the total transport cost is $\kappa = 15.47$. We then compute the intermediate $\mathcal{Y}_1$, midway along the homotopy path (i.e., 0.5 on the horizontal axis). This intermediate distribution is basically $\mathcal{Y}_1 = \mathcal{Y} {Q^*}^{1/2}$. The cost of using $T_2$ to transport $\mathcal{X}$ to $\mathcal{Y}_1$ is 33.63 (see the blue line at 0.5 on the homotopy path). $T_2$ is optimal to transport $\mathcal{X}$ to $\mathcal{Y}_2$ , not to transport $\mathcal{X}$ to $\mathcal{Y}_1$. Next step is to use $\mathcal{F}$ to optimize the transport plan between $\mathcal{X}$ and $\mathcal{Y}$, using $T_2$ as the initial solution. This yields $T_1$, corresponding to the drop in $\kappa$ from 33.63 to 20.86, at point 0.5 on the homotopy path. In the next homotopy iteration, we transform $\mathcal{Y}_1$ to $\mathcal{Y}$. Using $T_1$ as the transport plan between $\mathcal{X}$ and $\mathcal{Y}$ leads to $\kappa = 37.59$, as the blue line reaches the far right of horizontal axis. Using $\mathcal{F}$ to improve $T_1$ yields $T_0$ which is the optimal transport plan between $\mathcal{X}$ and $\mathcal{Y}$, leading to $\kappa = 29.06$. The improvement from 37.59 to 29.06 is depicted by the final drop in blue line at the far right of horizontal axis in Figure~\ref{fig:toy_pathsteps}, marking the end of this homotopy path.

When we make the discretization of the path finer, clearly, the solutions we trace are much closer and as a result, there are less fluctuations along the path. When $h$ is increased to $8$, the fluctuations are very small and the cost is changing almost linearly along the path, so we do not see a need to refine the discretization to $h = 16$.

\section{Other experiments}

We further experiment on a dataset of randomly generated data with $n=5,000$ and $d=500$, and observe that our homotopy algorithm is faster than the cubic running time algorithm, while it finds the same optimal solution. This trend repeated when we experiment on MNIST/CIFAR images, in the setting defined by~\cite{guo2020fast}.

\section{Conclusion and future work}

We presented the outline of a homotopy algorithm for the optimal transport problem with overall complexity of $\mathcal{O}(n^2 \log(n))$. To complete this work, we need to prove that the updates along the homotopy path have complexity no more than $\mathcal{O}(n^2 \log(n))$. This seems plausible, because there are algorithms in the literature that can provide updates in near-linear time with some approximation guarantees \citep{altschuler2017near}. We aim to improve this by proving that the updates along the homotopy path can be optimal with no approximation. The best known practical complexity bound for this problem is $\mathcal{O}(n^3)$ \cite{guo2020fast}.

\section*{Acknowledgements}
R.Y. was supported by a fellowship from the Department of Veteran Affairs. The views expressed in this manuscript are those of the author and do not necessarily reflect the position or policy of the Department of Veterans Affairs or the United States government.




\bibliography{refs}


\end{document}